\documentclass{article}
\usepackage{stmaryrd}
\usepackage{mathrsfs}
\usepackage[centertags]{amsmath}
\usepackage{amsfonts, dsfont}
\usepackage{amssymb}
\usepackage{amsthm}

\input xy
\xyoption{all}

\theoremstyle{plain}
\newtheorem{thm}{Theorem}[section]
\newtheorem{cor}[thm]{Corollary}
\newtheorem{lem}[thm]{Lemma}
\newtheorem{prop}[thm]{Proposition}

\theoremstyle{definition}
\newtheorem{defn}[thm]{Definition}
\newtheorem{remark}[thm]{Remark}
\newtheorem*{ack}{Acknowledgments}

\newcommand{\bd}{\begin{defn}}
\newcommand{\ed}{\end{defn}}
\newcommand{\bl}{\begin{lem}}
\newcommand{\el}{\end{lem}}
\newcommand{\bp}{\begin{prop}}
\newcommand{\ep}{\end{prop}}
\newcommand{\bt}{\begin{thm}}
\newcommand{\et}{\end{thm}}
\newcommand{\bc}{\begin{cor}}
\newcommand{\ec}{\end{cor}}
\newcommand{\br}{\begin{remark}}
\newcommand{\er}{\end{remark}}
\newcommand{\bdi}{\begin{diagram}}
\newcommand{\edi}{\end{diagram}}
\newcommand{\beq}{\begin{equation}}
\newcommand{\eeq}{\end{equation}}
\newcommand{\ba}{\begin{array}}
\newcommand{\ea}{\end{array}}
\newcommand{\bpf}{\begin{proof}}
\newcommand{\epf}{\end{proof}}

\newcommand{\Z}{\mathds{Z}}
\newcommand{\Q}{\mathds{Q}}
\newcommand{\Zp}{\mathds{Z}_{p}}
\newcommand{\Qp}{\mathds{Q}_{p}}
\newcommand{\al}{\alpha}
\newcommand{\Ga}{\Gamma}
\newcommand{\ga}{\gamma}
\newcommand{\La}{\Lambda}
\newcommand{\la}{\lambda}

\newcommand{\M}{\mathfrak{M}}

 \DeclareMathOperator{\Gal}{Gal}
\DeclareMathOperator{\Hom}{Hom} \DeclareMathOperator{\rank}{rank}

\DeclareMathOperator{\Ext}{Ext} \DeclareMathOperator{\Ann}{Ann}

\newcommand{\ot}{\otimes}
\newcommand{\ilim}{\displaystyle \mathop{\varinjlim}\limits}
\newcommand{\plim}{\displaystyle \mathop{\varprojlim}\limits}

\newcommand{\coker}{\mathrm{coker}\,}

\newcommand{\cyc}{\mathrm{cyc}}
\newcommand{\cts}{\mathrm{cts}}

\newcommand{\lra}{\longrightarrow}

\newcommand{\ps}[1]{\llbracket #1 \rrbracket}

\begin{document}

\title{On the complete faithfulness of the
   $p$-free quotient modules of dual Selmer groups}
\author{Meng Fai Lim\footnote{School of Mathematics and
Statistics, Central China Normal University, 152 Luoyu Road, Wuhan,
Hubei, P.R.China 430079. E-mail: \texttt{limmf@mail.ccnu.edu.cn}}}
\date{}
\maketitle

\begin{abstract} \footnotesize
\noindent
 In this paper, we consider the question of the
complete faithfulness of the $p$-free quotient module of the dual
Selmer groups of elliptic curves defined over a noncommutative
$p$-adic Lie extension. Our question will refine previous questions
on the complete faithfulness of dual Selmer groups. We also
consider the question of the triviality of the central torsion
submodules of these Iwasawa modules and we see that this latter question is intimately related to the former. We will also formulate and study
analogous questions for the dual Selmer groups of Hida deformations.
We then give positive answer to our questions, and establish ``control theorem" results between the questions in certain cases.



\end{abstract}

\section{Introduction}

Let $p$ be an odd prime. We are interested in a certain class of
modules defined over the Iwasawa algebra
 \[ \Zp\ps{G} = \plim_U \Zp[G/U] \]
of a noncommutative compact pro-$p$ torsionfree Lie group $G$. The
modules belonging to this class are called completely faithful
modules (see Section \ref{Algebraic Preliminaries} for the precise
definition). The motivation of studying these modules arises from
noncommutative Iwasawa theory in the search of a
global annihilator for the dual Selmer group of an elliptic curve
defined over a noncommutative $p$-adic Lie extension. Such a global
annihilator, if it exists, will give insight to the noncommutative
$p$-adic $L$-function. Unfortunately, the attempt to find a global annihilator proved futile, and this was first pointed out by Hachimori and
Venjakob in \cite{HV}. There, by building on a previous work of
Venjakob \cite{V03}, they constructed examples of completely
faithful dual Selmer groups over a certain class of
$\Zp\rtimes\Zp$-extensions. We should mention that the question of
the complete faithfulness of the dual Selmer group was also raised
in \cite{CSS}. Since then, many authors, including the present
author, have studied completely faithful modules over Iwasawa
algebra of other compact $p$-adic Lie groups and constructed
examples of completely faithful dual Selmer groups over $p$-adic Lie
extensions with Galois groups realized by such $p$-adic Lie groups
(see \cite{A, BZ, Cs, LimCF}). In this paper, we will consider a
refinement of the question of complete faithfulness which we
now describe.

Let $F$ be a number field and $E$ an elliptic curve
defined over $F$ which has either good ordinary reduction or
multiplicative reduction at each prime of $F$ above $p$. Let
$F_{\infty}$ be a $p$-adic Lie extension of $F$ with the properties
that the Galois group $G= \Gal(F_{\infty}/F)$ is a noncommutative
compact pro-$p$ Lie group with no $p$-torsion and that $F_{\infty}$
contains the cyclotomic $\Zp$-extension $F^{\cyc}$ of $F$. Write
$H=\Gal(F_{\infty}/F^{\cyc})$. A common phenomenon in all the
examples of completely faithful dual Selmer groups mentioned in the previous paragraph is that they are all finitely
generated over $\Zp\ps{H}$. In fact, the question raised in
\cite{CSS} was also posed under the condition when the dual Selmer
group in question is finitely generated over $\Zp\ps{H}$. It is then
natural to ask what can be said if the dual Selmer group is not
finitely generated over $\Zp\ps{H}$. Denote by $X(E/F_{\infty})$
the dual Selmer group of $E$ defined over $F_{\infty}/F$. We write
$X(E/F_{\infty})(p)$ for the $p$-primary submodule of $X(E/F_{\infty})$
and write $X_f(E/F_{\infty}) = X(E/F_{\infty})/X(E/F_{\infty})(p)$.
As we will see in Lemma \ref{mu lemma}(d), a necessary condition for
a $\Zp\ps{G}$-module $M$ to be completely faithful is that its $p$-primary
submodule must be pseudo-null. As is well-known, there are examples of
$X(E/F_{\infty})(p)$ being not pseudo-null. Building on
this observation, the author has produced examples of faithful dual Selmer
groups that are not completely faithful in \cite{LimCF}. Therefore, one may cut off the
$p$-primary submodule of $X(E/F_{\infty})$, and ask the natural
question of whether $X_f(E/F_{\infty})$ is completely faithful. This
will be the main theme of the paper. Throughout the paper,
$X_f(E/F_{\infty})$ is said to satisfy Property CF whenever
$X_f(E/F_{\infty})$ is either pseudo-null or completely faithful,
and we will be interested in the question of the validity of
Property CF.
We also propose and study a question on the triviality of the
central torsion submodules of the dual Selmer group, and we say that $X_f(E/F_{\infty})$ satisfies
Property CT whenever this triviality condition holds.
As already observed in literature (see \cite{A,BZ,CSS,Cs}), these two questions are intimately related.
We also formulate versions of Property CF and
Property CT for an appropriate Selmer group attached to an ordinary Hida deformation
(see Section \ref{completely faithful Hida section}). We mention that the variant of the Selmer group used here was first introduced by Greenberg \cite{G89} and has been incorporated to the noncommutative situation (for insatnce, see \cite{CS12,SS}). While formulating these properties, we also take the opportunity to collect and review results from existing literature to establish some instances of the validity of these properties.

We now describe the main results of this paper. As mentioned in the previous paragraph, Property CF and Property CT are intimately related. In most arithmetic situations, one can deduce the validity of Property CT from the validity of CF (see Sections \ref{triviality of central torsion submodules section} and \ref{completely faithful Hida section}). The reverse implication is a slightly more subtle issue, and there has been some recent studies of Ardakov \cite{A} and Csige \cite{Cs} in this direction. In this paper, we will establish a ``control theorem" type result
which allows the deduction of the validity of Property CT for an
elliptic curve from the validity of Property CF over a
\textit{smaller} extension (see Theorem \ref{faithful trivial
torsion}). We also establish an analogous result which allows us to deduce
the validity of Property CT for the Selmer group of the Hida deformation from the
validity of Property CF for the Selmer group of its specialization (see Theorem
\ref{faithful control for Hida CT}). The point of our result is that we do not have a control theorem for Property CF in general yet. In a previous paper \cite{LimCF}, the author established a control theorem for faithfulness of the Selmer group which is a weaker assertion than Property CF (see \cite[Section 6]{LimCF}). In view of the theorem obtained there, it is natural to believe (and hope) that there should be a control theorem for Property CF. Since our results are implied by such a control theorem, they can be viewed as partial support to this belief. Furthermore, by combining our ``control
theorems" with the results of Ardakov \cite{A} and Csige \cite{Cs}, we establish a control theorem for Property CF (see Theorems \ref{completely faithful control} and \ref{faithful control
for Hida}) for a certain class of $p$-adic Lie extensions, thus giving a first concrete substantiation of our belief.

We now give a brief description of the layout of the paper. In Section \ref{Algebraic Preliminaries},
we will review certain algebraic preliminaries which will be used subsequently in the paper. In Section
\ref{completely faithful section} (resp., Section \ref{triviality of central torsion submodules section}), we formulate Property CF (resp., Property CT) for the Selmer group of an elliptic curve. The analogue of Property CF and Property CT for the Selmer group of a Hida deformation will be introduced in Section \ref{completely faithful Hida section}. In Section \ref{control CT section}, we will prove our main ``control theorems". We then combine these theorems with the results of Ardakov and Csige to obtain control theorems for Property CF in Section \ref{control section}. Finally, in Section \ref{examples section}, we discuss some numerical examples.

\section{Algebraic Preliminaries} \label{Algebraic Preliminaries}

In this section, we establish some more or less standard algebraic
preliminaries and notation which are necessary for the discussion in the paper. Throughout
the paper, we will always work with left modules over a ring, where the ring is in turn always assumed to contain 1. Let $\La$ be a Auslander regular ring (for instance, see
\cite[Definition 3.3]{V02}) with no zero divisors. In particular, $\La$ is a (not
neccessarily commutative) Noetherian ring which has no zero
divisors. It is then well-known that the ring $\La$ admits a skew field of
fractions $K(\La)$ which is flat over $\La$ (see \cite[Chapters 6
and 10]{GW} or \cite[Chapter 4, \S 9 and \S 10]{Lam}). Thanks to this well-known fact, one can define the notion of $\La$-rank and $\La$-torsion module which we now do. For a
finitely generated $\La$-module $M$, its $\La$-rank is defined by the following formula:
$$ \rank_{\La}M  = \dim_{K(\La)} K(\La)\ot_{\La}M. $$
Then the $\La$-module $M$ is said to be \textit{torsion} if $\rank_{\La}M
=0$.

We now set to define the notion of a completely faithful $\La$-module. To prepare for this, we first recall that for a nonzero $\La$-module $M$, the global annihilator ideal of $M$
is defined by
\[ \Ann_{\La}(M) = \{ \la\in \La : \la m = 0~\mbox{for all}~m\in
M\}.\]
It is then a fairly straightforward exercise to check that this is a two-sided ideal of $\La$, and we
will say that $M$ is a \textit{faithful} $\La$-module if
$\Ann_{\La}(M) = 0$. Denote by $\mathcal{M}$ the category of
all finitely generated torsion $\La$-modules and by $\mathcal{C}$
the full subcategory of all pseudo-null $\La$-modules in $\mathcal{M}$. Here we recall that a
finitely generated torsion $\La$-module $M$ is said to be \textit{pseudo-null}
if $\Ext^1_{\La}(M,\La)=0$.  Write $q :\mathcal{M}\lra \mathcal{M}/\mathcal{C}$ for the quotient
functor.  Equipped with these notions, a finitely generated torsion $\La$-module $M$ is then said to be \textit{completely faithful} if $\Ann_{\La}(N) = 0$ for
any $N \in\mathcal{M}$ such that $q(N)$ is isomorphic to a non-zero
subquotient of $q(M)$. The following simple observation will be frequently used in the paper.

\bl \label{cf compare} Suppose that
\[ \al : M\lra N\]
is a $\La$-homomorphism of finitely generated torsion $\La$-modules
with the property that $\ker\al$ and $\coker\al$ are pseudo-null
$\La$-modules, then $M$ is a completely faithful $\La$-module
$($resp., pseudo-null $\La$-module$)$ if and only if $N$ is a
completely faithful $\La$-module $($resp., pseudo-null
$\La$-module$)$. \el

\bpf
 Since $\ker\al$ and $\coker\al$ are pseudo-null
$\La$-modules, we have that $q(\al):q(M)\lra q(N)$ is an isomorphism in $\mathcal{M}/\mathcal{C}$.
The assertion of the lemma is now immediate from this.
\epf

We introduce further notation which we will frequently use in this paper without
further comment. Let $x\in \La$. Denote by $M[x]$ the set consisting of
elements of $M$ annihilated by $x$. In general, $M[x]$ is
at most an additive subgroup of $M$. However, if one assumes further
that $x\La = \La x$, it is then easy to verify that
$M[x]$ is a $\La$-submodule of $M$. Also, one can check
easily that
 \[xM = \{xm : m\in M\}  \]
is a $\La$-submodule of $M$, and that $M/xM$ is a $\La/x\La$-module.

We now introduce the class of Auslander regular rings with no zero divisors that
will be considered in this paper. Namely, this class comes from the Iwasawa algebras of
compact $p$-adic Lie groups which was first discovered by Venjakob \cite{V02}. From now on, $p$ will denote a fixed odd prime, and
$G$ will denote a compact pro-$p$ $p$-adic Lie group without $p$-torsion. The
completed group algebra of $G$ over $\Zp$ is given by
 \[ \Zp\ps{G} = \plim_U \Zp[G/U], \]
where $U$ runs over the open normal subgroups of $G$ and the inverse
limit is taken with respect to the canonical projection maps.
A celebrated theorem of Lazard \cite{Laz} asserts that the ring $\Zp\ps{G}$ is Noetherian.
Venjakob went one step further by establishing that $\Zp\ps{G}$ is an
Auslander regular ring (see \cite[Theorem 3.26]{V02}). Furthermore, it follows from an earlier work of Neumann \cite{Neu} that
the ring $\Zp\ps{G}$ has no zero divisors. Therefore, there is a well-defined notion of
$\Zp\ps{G}$-rank, torsion $\Zp\ps{G}$-module, pseudo-null
$\Zp\ps{G}$-module and completely faithful $\Zp\ps{G}$-module.
We record another well-known and important result of
Venjakob (cf. \cite[Example 2.3 and Proposition 5.4]{V03}) which we
require in the discussion of the paper.

\bt[Venjakob] \label{pseudo-null torsion} Suppose that $H$ is a
closed normal subgroup of $G$ with $G/H \cong \Zp$. For a
compact $\Zp\ps{G}$-module $M$ which is finitely generated over
$\Zp\ps{H}$, we have that $M$ is a pseudo-null $\Zp\ps{G}$-module if and
only if $M$ is a torsion $\Zp\ps{H}$-module. \et

We record another useful lemma which is a special case of \cite[Lemma 4.5]{LimFine}.

\bl \label{relative rank} Let $H$ be a compact pro-$p$ $p$-adic Lie
group without $p$-torsion. Let $N$ be a closed normal subgroup of
$H$ such that $N\cong \Zp$ and such that $H/N$ is also a compact
pro-$p$ $p$-adic Lie group without $p$-torsion. Let $M$ be a
finitely generated $\Zp\ps{H}$-module. Then $H_i(N,M)$ is finitely
generated over $\Zp\ps{H/N}$ for each $i$ and $H_i(N,M) =0$ for
$i\geq 2$. Furthermore, we have an equality
\[\rank_{\Zp\ps{H}}M = \rank_{\Zp\ps{H/N}}M_N - \rank_{\Zp\ps{H/N}}H_1(N,M). \] \el

Now, for a given finitely generated $\Zp\ps{G}$-module $M$, denote
by $M(p)$ the $\Zp\ps{G}$-submodule of $M$ consisting of elements of
$M$ which are annihilated by some power of $p$. We shall call $M(p)$
the $p$-primary submodule of $M$. A finitely generated
$\Zp\ps{G}$-module $M$ is then said to be \textit{$p$-primary} if
$M= M(p)$. Since the ring $\Zp\ps{G}$ is Noetherian and $M$ is
finitely generated over $\Zp\ps{G}$, the module $M(p)$ is certainly
finitely generated over $\Zp\ps{G}$. As a consequence, there exists
a sufficiently large integer $r\geq 0$ such that $p^r$ annihilates
$M(p)$. Following \cite[Formula (33)]{Ho}, we define
  \[\mu_G(M) = \sum_{i\geq 0}\rank_{\mathbb{F}_p\ps{G}}\big(p^i
   M(p)/p^{i+1}\big), \]
which is a finite sum by our discussion. (For another alternative, but equivalent, definition, see
\cite[Definition 3.32]{V02}.) It is immediate from the definition that $\mu_{G}(M) =
\mu_{G}(M(p))$.
The following lemma collects several properties of the $\mu_G$-invariant which
will be used in the paper.

\bl \label{mu lemma} Let $G$ be a compact pro-$p$ $p$-adic Lie group
without $p$-torsion and $M$ a finitely generated $\Zp\ps{G}$-module.
Then we have the following statements.

\begin{enumerate}
\item[$(a)$] The $\mu_G$-invariant is additive on short exact
sequences of torsion $\Zp\ps{G}$-modules.

\item[$(b)$]  Suppose that $G$ has a closed normal subgroup $H$ such that
$G/H\cong \Zp$. If $M$ is a $\Zp\ps{G}$-module which is finitely
generated over $\Zp\ps{H}$, then $\mu_G(M) =0$.

\item[$(c)$] We have $\mu_G(M)
=0$ if and only if $M(p)$ is pseudo-null over $\Zp\ps{G}$. In
particular,  if $M$ is a $p$-primary $\Zp\ps{G}$-module, then
$\mu_G(M) =0$ if and only if $M$ is pseudo-null over $\Zp\ps{G}$.

\item[$(d)$] Suppose that $M$ is a completely faithful $\Zp\ps{G}$-module.
Then $M(p)$ is pseudo-null over $\Zp\ps{G}$ and $\mu_G(M) = 0$.
 \end{enumerate}
\el

\bpf Statement (a) follows from \cite[Proposition 1.8]{Ho} (see also
\cite[Corollary 3.37]{V02}). Statement (b) is precisely \cite[Lemma
2.7]{Ho}. The validity of statement (c) is established in \cite[Remark 3.33]{V02}. We now
show statement (d). Since $M(p)$ is a $\Zp\ps{G}$-submodule of $M$
which is annihilated by some power of $p$ and $M$ is completely
faithful, we must have that $M(p)$ is a pseudo-null $\Zp\ps{G}$-module.
 By statement (c), this in turn implies that $\mu_G(M) =
 \mu_G(M(p)) = 0$.
  \epf

\section{A question on complete faithfulness} \label{completely
faithful section}

In this section, we introduce and formulate Property CF for the Selmer group of an elliptic curve.
Fix once and for all an algebraic closure $\bar{\Q}$ of $\Q$.
Therefore, an algebraic (possibly infinite) extension of $\Q$ will
mean an subfield of $\bar{\Q}$. A finite extension $F$ of $\Q$ will
be called a number field. Let $E$ be an elliptic curve which is
defined over a number field $F$. Suppose that for every prime $v$ of
$F$ above $p$, the elliptic curve $E$ has either good ordinary
reduction or multiplicative reduction at $v$.

Let $v$ be a prime of $F$. For every finite extension $L$ of $F$, we
define
 \[ J_v(E/L) = \bigoplus_{w|v}H^1(L_w, E)_{p^{\infty}},\]
where $w$ runs over the (finite) set of primes of $L$ above $v$. For an infinite algebraic extension $\mathcal{L}$ of $F$, set
\[ J_v(E/\mathcal{L}) = \ilim_L J_v(E/L),\]
where the direct limit is taken over all finite extensions $L$ of
$F$ contained in $\mathcal{L}$. The Selmer group of $E$
over $\mathcal{L}$ is then denoted and defined by
\[ S(E/\mathcal{L}) = \ker\Big(H^1(\mathcal{L}, E_{p^{\infty}})
\lra \bigoplus_{v} J_v(E/\mathcal{L}) \Big), \] where $v$ runs
through all the primes of $F$.

In this paper, we are interested in studying Selmer groups defined over an important
class of infinite algebraic extensions of $F$ which we now define.
A Galois extension $F_{\infty}$ of $F$ is said to be an \textit{admissible
$p$-adic Lie extension} of $F$ if (i) $\Gal(F_{\infty}/F)$ is a
compact $p$-adic Lie group, (ii) $F_{\infty}$ contains the
cyclotomic $\Zp$-extension $F^{\cyc}$ of $F$ and (iii) $F_{\infty}$
is unramified outside a finite set of primes of $F$. Furthermore, if the group
$\Gal(F_{\infty}/F)$ has no $p$-torsion, we say that $F_{\infty}$ is a \textit{strongly
admissible $p$-adic Lie extension} of $F$.  We shall write $G =
\Gal(F_{\infty}/F)$, $H = \Gal(F_{\infty}/F^{\cyc})$ and $\Ga
=\Gal(F^{\cyc}/F)$. Let $S$ be any finite set
of primes of $F$ which contains the primes above $p$, the infinite
primes, the primes at which $E$ has bad reduction and the primes
that are ramified in $F_{\infty}/F$. The maximal
algebraic extension of $F$ unramified outside $S$ is then denoted by $F_S$, and for every
algebraic extension $\mathcal{L}$ of $F$
contained in $F_S$, we write $G_S(\mathcal{L}) =
\Gal(F_S/\mathcal{L})$. The following equivalent
description 
\[ S(E/F_{\infty}) = \ker\Big(H^1(G_S(F_{\infty}), E_{p^{\infty}})
\stackrel{\lambda_{E/F_{\infty}}}{\lra} \bigoplus_{v\in S}
J_v(E/F_{\infty}) \Big)\] of the Selmer group of $E$ over $F_{\infty}$ is well-known (for instance, see
\cite[Lemma 2.2]{CS12}) and will be utilized without further comment. From now on, we shall
 denote by $X(E/F_{\infty})$ the
Pontryagin dual of $S(E/F_{\infty})$.

As mentioned in the Introduction, when $X(E/F_{\infty})$ is
finitely generated over $\Zp\ps{H}$, the question of the complete faithfulness of $X(E/F_{\infty})$ 
has been studied by many authors (see \cite{A,BZ,CSS,HV,LimCF}). It is then natural to
consider the situation $X(E/F_{\infty})$ is not finitely
generated over $\Zp\ps{H}$. In this case, the present author has previously
constructed examples of dual Selmer groups which are faithful but
not completely faithful (see \cite{LimCF}). As observed
loc.\ cit., these examples of dual Selmer groups
have positive $\mu$-invariants. Furthermore, in view of Lemma \ref{mu
lemma}(d), to be able to speak of complete
faithfulness of a general module, one needs to cut off its $p$-primary submodule. We
will therefore introduce the following property which will be the
main theme of the paper.

Write $X_f(E/F_{\infty}) :=
X(E/F_{\infty})/X(E/F_{\infty})(p)$. We say that $X(E/F_{\infty})$ (or $X_f(E/F_{\infty})$) \textit{satisfies Property CF} if
$X_f(E/F_{\infty})$ is either a pseudo-null
$\Zp\ps{G}$-module or a completely faithful $\Zp\ps{G}$-module.
In this paper, we will be
interested in studying the validity of Property CF. We first show that this question covers the
original question of Coates et al.

\bl \label{Coates}
 Let $E$ be an elliptic curve over $F$
which has either good ordinary reduction or multiplicative reduction
at every prime of $F$ above $p$.
 Let $F_{\infty}$ be a noncommutative strongly admissible $p$-adic Lie extension of $F$ with
 $G =\Gal(F_{\infty}/F)$. Assume that $X(E/F_{\infty})$ is finitely generated over
  $\Zp\ps{H}$. Then the following statements are equivalent.
 \begin{enumerate}
 \item[$(a)$] $X(E/F_{\infty})$ is either a pseudo-null
$\Zp\ps{G}$-module or a completely faithful $\Zp\ps{G}$-module.
  \item[$(b)$] $X_f(E/F_{\infty})$ is either a pseudo-null
$\Zp\ps{G}$-module or a completely faithful $\Zp\ps{G}$-module.
 \end{enumerate}
\el

\bpf
 Since $X(E/F_{\infty})$ is finitely generated over
  $\Zp\ps{H}$, so is $X(E/F_{\infty})(p)$. It then follows from Lemma \ref{mu lemma}(b) and (c) that
  $X(E/F_{\infty})(p)$ is pseudo-null over $\Zp\ps{G}$. The
  equivalence of the statements of the lemma is now
  an immediate consequence of Lemma \ref{cf
  compare}.
 \epf

\br \label{00}
 In fact, if $E$ has good ordinary reduction at all the primes
 above $p$ and $F_{\infty}$ is not totally real, we have the following sharper observation. Namely, if $X(E/F_{\infty})$ is
 finitely generated over $\Zp\ps{H}$, then $X(E/F_{\infty})(p) =
 0$ and $X(E/F_{\infty})= X_f(E/F_{\infty})$. To see this, first note that since $X(E/F_{\infty})$ is finitely generated over
$\Zp\ps{H}$, we have that $X(E/L^{\cyc})$ is finitely generated
over $\Zp$ for every finite extension $L$ of $F$ contained in
$F_{\infty}$. Since $F_{\infty}$ is not totally real, we can rewrite
$X(E/F_{\infty}) = \plim_LX(E/L^{\cyc})$, where $L$ runs through all
finite non totally real extensions of $F$ contained in $F_{\infty}$. Since $X(A/L^{\cyc})$ is finitely generated
over $\Zp$, by an application of \cite[Proposition
 7.5]{Mat}, the multiplication by $p$-map on $X(E/L^{\cyc})$ is an
injection. As inverse limit is left exact, it then follows that the
multiplication by $p$-map on $X(E/F_{\infty})$ is also an injection.
Consequently, we have $X(E/F_{\infty})(p) =
 0$ and $X(E/F_{\infty})= X_f(E/F_{\infty})$.
\er

The next proposition shows that the validity of Property CF is invariant under isogeny.

\bp \label{faithful isogeny}
 Let $E_1$  and $E_2$ be two elliptic curves over $F$
with either good ordinary reduction or multiplicative reduction at
every prime of $F$ above $p$. Let
$F_{\infty}$ be a noncommutative strongly admissible noncommutative
$p$-adic Lie extension of $F$ with $G= \Gal(F_{\infty}/F)$.
Suppose that $E_1$ and $E_2$ are isogenous to each other over $F$, and suppose that $X(E_1/F_{\infty})$ is a torsion $\Zp\ps{G}$-module. Then
$X(E_2/F_{\infty})$ is a torsion $\Zp\ps{G}$-module. Furthermore,
$X(E_1/F_{\infty})$ satisfies Property CF if and only if $X(E_2/F_{\infty})$ satisfies Property CF. \ep

\bpf Let $\varphi: E_2\lra E_1$ be an isogeny defined over $F$. Then $\varphi$ induces homomorphisms on the cohomology groups and Selmer groups, giving rise to the following commutative diagram 
\[
 \entrymodifiers={!! <0pt, .8ex>+} \SelectTips{eu}{}\xymatrix{
    0 \ar[r]^{} & S(E_2/F_{\infty}) \ar[d]_{\varphi_1} \ar[r] &
    H^1(G_S(F_{\infty}), E_{2,p^{\infty}})
    \ar[d]_{\varphi_2}
    \ar[r] & \bigoplus_{v\in S}J_v(E_2/F_{\infty}) \ar[d]_{\varphi_3} \\
    0 \ar[r]^{} & S(E_1/F_{\infty})\ar[r]^{}
    & H^1(G_S(F_{\infty}),  E_{1,p^{\infty}}) \ar[r] & \
    \bigoplus_{v\in S}J_v(E_1/F_{\infty})  } \]
    with exact rows.
Since $\ker\varphi_2 = H^1(G_S(F_{\infty}),\ker\varphi)$ and $\coker\varphi_2\subseteq  H^2(G_S(F_{\infty}),\ker\varphi)$, it follows that $\ker\varphi_2$ and $\coker\varphi_2$ are killed by the
$p$-order of $\ker\varphi$. Similarly, one can also show that $\ker\varphi_3$ is killed by the
$p$-order of $\ker\varphi$. Hence, by the snake lemma, $\ker\varphi_1$ is killed by the $p$-order of $\ker\varphi$, and $\coker\varphi_1$ is killed by a square of the $p$-order of $\ker\varphi$. Now consider the following commutative diagram
 \[ \entrymodifiers={!! <0pt, .8ex>+}
\SelectTips{eu}{}\xymatrix{
    0 \ar[r]^{} & X(E_1/F_{\infty})(p) \ar[d]_f \ar[r] & X(E_1/F_{\infty})
    \ar[d]_g
    \ar[r] & X_f(E_1/F_{\infty}) \ar[d]_{h} \ar[r] & 0 \\
    0 \ar[r]^{} & X(E_2/F_{\infty})(p) \ar[r]^{} & X(E_2/F_{\infty})\ar[r] & \
    X_f(E_2/F_{\infty}) \ar[r] &  0  } \]
with exact rows, where the map $g$ is precisely the Pontryagin dual of $\varphi_1$. It then follows from the above discussion that $g$ has kernel and cokernel that are $p$-primary $\Zp\ps{G}$-modules. Since $X(E_1/F_{\infty})$ is a torsion
$\Zp\ps{G}$-module by the hypothesis, it
follows that $X(E_2/F_{\infty})$ is also a torsion
$\Zp\ps{G}$-module.

Now, clearly, the map $f$ has kernel and cokernel that are $p$-primary
$\Zp\ps{G}$-modules. Combining this with the above observation that the map $g$ has kernel and cokernel that are $p$-primary $\Zp\ps{G}$-modules, we conclude that the map $h$ also has
kernel and cokernel that are $p$-primary $\Zp\ps{G}$-modules. Since
$ X_f(E_1/F_{\infty})$ has no $p$-torsion by definition, this implies that
$\ker h=0$. Consequently, we have $\ker f = \ker g$ and
the following short exact sequence
 \[ 0\lra \coker f \lra \coker g \lra
\coker h \lra 0.\] On the other hand, it follows from the exact
sequence of torsion $\Zp\ps{G}$-modules
\[ 0 \lra \ker g \lra X(E_1/F_{\infty}) \lra X(E_2/F_{\infty}) \lra
\coker g \lra 0\] and Lemma \ref{mu lemma}(a) that
 \[ \mu_G(\ker g) - \mu_G(\coker g) =   \mu_G\big(X(E_1/F_{\infty})\big)- \mu_G\big(X(E_2/F_{\infty})\big).\]
 Similarly, we have
 \[ \mu_G(\ker f) - \mu_G(\coker f) =   \mu_G\big(X(E_1/F_{\infty})(p)\big)- \mu_G\big(X(E_2/F_{\infty})(p)\big).\]
 Since $\mu_G(M)=\mu_G\big(M(p)\big)$, the above two equalities combine to give
 \[ \mu_G(\ker g) - \mu_G(\coker g) = \mu_G(\ker f) - \mu_G(\coker f).\]
As we have already observed above that $\ker f = \ker g$, it then
follows that $\mu_G(\coker f)= \mu_G(\coker g)$ which in turn
implies that $\mu_G(\coker h) = 0$. But since $\coker h$ is $p$-primary,
we have that $\coker h$ is a pseudo-null $\Zp\ps{G}$-module by an application of Lemma \ref{mu lemma}(c). Therefore, we have shown that the map
\[ h: X_f(E_1/F_{\infty})\lra X_f(E_2/F_{\infty})\]
is injective with a pseudo-null cokernel. The final assertion of the
proposition now follows from an application of Lemma \ref{cf
  compare}. \epf

\br
 If one assumes that the $\M_H(G)$-Conjecture (see below) holds for
 $X(E_1/F_{\infty})$ and $X(E_2/F_{\infty})$, then one can give a
 shorter proof of the preceding proposition. However, we thought it may be of
 interest to be able to give a proof of the proposition without assuming the
 $\M_H(G)$-Conjecture. We also mention that it is still an
 open question to whether the $\M_H(G)$-Conjecture is invariant
 under isogeny (but see \cite[Proposition 5.6]{CFKSV} for a
 partial result in this direction).
\er

To facilitate further discussion, we need to introduce a conjecture on our Selmer group which was first proposed in \cite{CFKSV}. We say that the $\M_H(G)$-Conjecture holds for $X(E/F_{\infty})$ if $X_f(E/F_{\infty})$ is a finitely generated
$\Zp\ps{H}$-module. If one assumes the $\M_H(G)$-Conjecture and takes Theorem \ref{pseudo-null torsion}
into account, then the validity of Property CF is equivalent to saying that
$X_f(E/F_{\infty})$ is either a torsion $\Zp\ps{H}$-module or a
completely faithful $\Zp\ps{G}$-module. The following lemma is not
used in the paper but we have decided to include it, as we believe
that it is an interesting observation.

\bl \label{observe}
 Let $E$ be an elliptic curve over $F$
which has either good ordinary reduction or multiplicative reduction
at every prime of $F$ above $p$.
 Let $F_{\infty}$ be a noncommutative strongly admissible $p$-adic Lie extension of $F$ with
 $G =\Gal(F_{\infty}/F)$. Suppose that $X(E/F_{\infty})$ satisfies the $\M_H(G)$-Conjecture and that $X(E/F_{\infty})$
 is a completely faithful $\Zp\ps{G}$-module. Then $X(E/F_{\infty})$
 is finitely generated over $\Zp\ps{H}$.
\el

\bpf
  As $X(E/F_{\infty})$
 is assumed to be completely faithful, we may apply Lemma \ref{mu lemma}(d) to conclude that $\mu_G\big(X(E/F_{\infty})\big) = 0$. Since $X(E/F_{\infty})$ also satisfies the
 $\M_H(G)$-Conjecture, it follows from \cite[Proposition
 2.5]{CS12} that $X(E/L^{\cyc})$ is torsion over $\Zp\ps{\Ga_L}$ for every finite
 extension $L$ of $F$ contained in $F_{\infty}$, where $\Ga_L=\Gal(L^{\cyc}/L)$. By \cite[Corollary
 3.4]{LimMHG}, this in turn implies that all the assumptions of \cite[Theorem
 3.1]{LimMHG} are satisfied, and therefore, we may apply the said
 theorem to conclude that $\mu_{\Ga}\big(X(E/F^{\cyc})\big) = 0$.
 Combining this with the above observation of $X(E/F^{\cyc})$ being torsion over $\Zp\ps{\Ga}$,
 we have that $X(E/F^{\cyc})$ is finitely generated
 over $\Zp$ by the structure theory of
$\Zp\ps{\Ga}$-modules. The required conclusion now follows from an
application of \cite[Theorem 2.1]{CS12}. \epf

\br
 The conclusion of Lemma \ref{observe} is false if we replace
 $X(E/F_{\infty})$ by an arbitrary $\Zp\ps{G}$-module. Suppose that $G$ satisfies $\mathbf{(NH)}$ (see below).
 We shall now construct a class of modules which are
 completely faithful over $\Zp\ps{G}$ and are not finitely generated over $\Zp\ps{H}$.
 Let $M_0$ be any $\Zp\ps{G}$-module which is finitely generated over $\Zp\ps{H}$
 with positive $\Zp\ps{H}$-rank. Let $M$ be the module
 $M_0\oplus (\Zp\ps{\Ga}/p)$, where $G$ acts on $\Ga$ via the quotient
 $G\twoheadrightarrow \Ga$. Clearly, $M/M(p)$ is finitely generated
 over $\Zp\ps{H}$, and we have a short exact sequence
 \[ 0\lra M_0\lra M \lra \Zp\ps{\Ga}/p \lra 0. \]
 Clearly, $\Zp\ps{\Ga}/p$ has trivial $\mu_G$-invariant, and therefore, is a pseudo-null
 $\Zp\ps{G}$-module by Lemma \ref{mu lemma}(d). It then follows from Lemma \ref{cf compare} and
 Theorem \ref{completely faithful modules} that $M$ is completely
 faithful over $\Zp\ps{G}$. On the other hand, since the module $\Zp\ps{\Ga}/p$ is not finitely
 generated over $\Zp$ and $H$ acts trivially on $\Zp\ps{\Ga}/p$, it is not finitely generated over $\Zp\ps{H}$. Hence $M$ is also not finitely generated over $\Zp\ps{H}$.
\er

For the remainder of the section, we will show that the results in
\cite{LimCF,V03} can be applied to verify Property CF in certain
cases. We shall first recall the result that we require. The group $G$
is said to satisfy $\mathbf{(NH)}$ if it contains two closed
normal subgroups $N$ and $H$ with the following two properties:
  \begin{enumerate}
\item[$(i)$] $N\subseteq H$, $G/H\cong \Zp$ and
 $G/N$ is a non-abelian group isomorphic to $\Zp\rtimes\Zp$.
\item[$(ii)$] There is a finite family of closed normal subgroups $N_i$
of $G$ such that $1=N_0\subseteq N_1 \subseteq
\cdots\subseteq N_r =N$ and such that for
$1\leq i\leq r$, one has $N_i/N_{i-1}\cong \Zp$.
\end{enumerate}

\bp \label{completely faithful modules}
  Suppose that $G$ satisfies $\mathbf{(NH)}$.
 Let $M$ be a $\Zp\ps{G}$-module which
is finitely generated over $\Zp\ps{H}$ with positive
$\Zp\ps{H}$-rank. Then $M$ is a completely faithful
$\Zp\ps{G}$-module. \ep

\bpf
 When $r=0$, this is precisely \cite[Corollary
4.3]{V03}. For $r>0$, one proceeds by an inductive argument. We
refer readers to \cite[Theorem 3.3]{LimCF} for the details. \epf

We can now verify the validity of Property CF
when the strongly admissible extension has a Galois group of the
form $\mathbf{(NH)}$.

\bp \label{faithful Selmer groups}
 Let $E$ be an elliptic curve over $F$
which has either good ordinary reduction or multiplicative reduction
at every prime of $F$ above $p$.
 Let $F_{\infty}$ be a strongly admissible $p$-adic Lie extension of $F$ with
 $G =\Gal(F_{\infty}/F)$. Suppose that the following statements hold.
 \begin{enumerate}
 \item[$(i)$] The $\M_H(G)$ conjecture holds for $X(E/F_{\infty})$.
\item[$(ii)$] $G$ satisfies $\mathbf{(NH)}$.
\end{enumerate}
  Then $X_f(E/F_{\infty})$ satisfies Property CF. \ep

\bpf
 If $X_f(E/F_{\infty})$ is a torsion $\Zp\ps{H}$-module, then the conclusion of the proposition will follow from
 by an application
 of Theorem \ref{pseudo-null torsion}. In the case that $X_f(E/F_{\infty})$ is a finitely generated $\Zp\ps{H}$-module
 with positive $\Zp\ps{H}$-rank, we can apply Proposition \ref{completely faithful
 modules} to conclude
 that $X_f(E/F_{\infty})$ is a completely faithful $\Zp\ps{G}$-module.
\epf

At this point of writing, we feel that there is currently
insufficient evidence for a positive answer to the question of Property CF being
satisfied in general. Hence it may be premature to make a conjecture
on Property CF, and therefore, we have refrained from doing so.

\section{A question on the triviality of central torsion submodule} \label{triviality of central torsion submodules section}

Retain the notation of the Section \ref{completely faithful
section}. We now introduce another property of the dual Selmer
groups which we like to investigate in this paper. Denote by $C$ the center of $G = \Gal(F_{\infty}/F)$, where $F_{\infty}$ is a noncommutative strongly admissible $p$-adic Lie extension of $F$. We say that $X(E/F_{\infty})$ (or $X_f(E/F_{\infty})$) \textit{satisfies Property CT} if the $\Zp\ps{C}$-torsion submodule of $X_f(E/F_{\infty})$ is zero. We note that it has been shown in \cite[Theorem 6.5]{OcV02} that the
$\Zp\ps{C}$-torsion submodule $X(E/F_{\infty})$ is either trivial or
not finitely generated over $\Zp\ps{C}$.
The next proposition records the isogeny invariance of Property CT.

\bp \label{trivial torsion isogeny}
 Let $E_1$  and $E_2$ be two elliptic curves over $F$
with either good ordinary reduction or multiplicative reduction at
every prime of $F$ above $p$ which are isogenous to each other. Let
$F_{\infty}$ be a noncommutative strongly admissible noncommutative
$p$-adic Lie extension of $F$ with $G= \Gal(F_{\infty}/F)$. Assume
that $X(E_1/F_{\infty})$ is a torsion $\Zp\ps{G}$-module $($and
hence so is $X(E_2/F_{\infty})$$)$. Then
$X(E_1/F_{\infty})$ satisfies Property CT if and only if
$X(E_2/F_{\infty})$ satisfies Property CT. \ep

\bpf
 Let $\varphi: E_2\lra E_1$ be an isogeny defined over $F$. From the
 proof of Proposition \ref{faithful isogeny}, we see that $\varphi$
 induces an injection \[ h: X_f(E_1/F_{\infty})\lra X_f(E_2/F_{\infty}).\]
 Therefore, if the central torsion submodule of
 $X_f(E_2/F_{\infty})$ is trivial, so is the central torsion submodule of
 $X_f(E_1/F_{\infty})$. By considering the dual isogeny $\hat{\varphi}: E_1\lra
 E_2$, we obtain the reverse implication.
\epf

 We now mention how we can deduce the validity of
Property CT in certain situations by appealing to Proposition \ref{faithful Selmer groups}. As
a start, we begin with a general observation.

\bp \label{CF CT}
 Let $E$ be an elliptic curve over $F$
which has either good ordinary reduction or multiplicative reduction at
every prime of $F$ above $p$. Let $F_{\infty}$ be a noncommutative strongly admissible $p$-adic Lie extension of $F$. Suppose that the following statements hold.
 \begin{enumerate}
  \item[$(i)$] The $\M_H(G)$ conjecture holds for
 $X(E/F_{\infty})$.
 \item[$(ii)$] $X(E/F_{\infty})$ has no nonzero pseudo-null $\Zp\ps{G}$-submodule.
\item[$(iii)$] Suppose that either $(a)$ or $(b)$ holds.
 \begin{enumerate}
  \item[$(a)$] $X_f(E/F_{\infty})$ satisfies Property CF.
 \item[$(b)$] The center $C$ of $G=\Gal(F_{\infty}/F)$ is contained in
 $H$.
\end{enumerate}
\end{enumerate}
 Then $X_f(E/F_{\infty})$ satisfies Property CT.
 \ep

 \bpf
   Since
 $X(E/F_{\infty})$ has no nonzero pseudo-null $\Zp\ps{G}$-submodule,
 it follows from \cite[Lemma 4.2]{Su} that $X_f(E/F_{\infty})$ also has no
 nonzero pseudo-null $\Zp\ps{G}$-submodule. We first consider the case when condition (iii)(a) holds. If $X_f(E/F_{\infty})$ is
 pseudo-null, then it follows that $X_f(E/F_{\infty})=0$, and in particular,
its $\Zp\ps{C}$-torsion submodule also vanishes. Now suppose that
$X_f(E/F_{\infty})$ is completely faithful over $\Zp\ps{G}$. As
$X_f(E/F_{\infty})$ is a Noetherian $\Zp\ps{G}$-module, so is its
$\Zp\ps{C}$-torsion submodule. As a consequence, the
$\Zp\ps{C}$-torsion submodule of $X_f(E/F_{\infty})$ has a global
annihilator, and therefore, must be pseudo-null over $\Zp\ps{G}$ by
the complete faithfulness of $X_f(E/F_{\infty})$. But since we have shown that
$X_f(E/F_{\infty})$ has no nonzero pseudo-null
$\Zp\ps{G}$-submodule, this in turn implies that the $\Zp\ps{C}$-torsion
submodule of $X_f(E/F_{\infty})$ is zero.

We now consider the case when condition (iii)(b) holds. Let $W$ be
the $\Zp\ps{C}$-torsion submodule of $X_f(E/F_{\infty})$. Clearly,
$W$ is a $\Zp\ps{G}$-submodule, and at the same time, a
$\Zp\ps{H}$-submodule of $X_f(E/F_{\infty})$. Since
$X_f(E/F_{\infty})$ is finitely generated over $\Zp\ps{H}$, we can
find $z\in \Zp\ps{C}$ such that $zW=0$. Since $\Zp\ps{C}\subseteq
Z\ps{H}$ by condition (iii)(b), this in turn implies that $W$ is a torsion $\Zp\ps{H}$-module. Hence, by
Theorem \ref{pseudo-null torsion}, $W$ is a pseudo-null
$\Zp\ps{G}$-module. But we have shown above that
$X_f(E/F_{\infty})$ has no nonzero pseudo-null
$\Zp\ps{G}$-submodule, and therefore, we must have $W=0$. \epf

We mention that condition (ii) of the preceding proposition is known to hold
in many cases when $E$ has good ordinary reduction at
every prime of $F$ above $p$ (see \cite[Theorem 3.2]{HO},
\cite[Theorem 2.6]{HV} and \cite[Theorem 5.1]{OcV02}). Hence we can combine these results with Theorem \ref{faithful
Selmer groups} to obtain many cases where Property CT is satisfied.

In the event when we do not know the validity of condition (ii), we
have the following result which will be useful in the discussion of
most numerical examples.

\bp \label{CF CT2}
 Let $E$ be an elliptic curve over $F$
which has good ordinary reduction at every prime of $F$ above $p$.
 Let $F_{\infty}$ be a noncommutative strongly admissible $p$-adic Lie extension
 of $F$ which is not totally real. Assume that the center $C$ of $G=\Gal(F_{\infty}/F)$ is contained in
 $H$.
 Suppose that there is an elliptic curve $E'$ which is
 isogenous to $E$ over $F$ and has the property that $X(E'/F^{\cyc})$ is
 finitely generated over $\Zp$.
 Then $X_f(E/F_{\infty})$ satisfies Property CT . \ep

 \bpf
   By virtue of
 Proposition \ref{trivial torsion isogeny} and Remark \ref{00}, it suffices to show
 that the $\Zp\ps{C}$-torsion submodule of $X(E'/F_{\infty})$ is zero. Also, it follows from
 \cite[Remark 2.6 and Corollary 6.1 (on p.
 1243)]{BZ} that $X(E'/F_{\infty})$ has no nonzero pseudo-null
$\Zp\ps{G}$-submodule. The conclusion now follows from a similar
argument to that in Proposition \ref{CF CT}.
 \epf

\br Note that since the property of having no
nonzero pseudo-null $\Zp\ps{G}$-submodule is not known to be isogeny invariant, we
cannot deduce that $X(E/F_{\infty})$ has no nonzero pseudo-null
$\Zp\ps{G}$-submodule. \er

\section{Analogue of Property CF and Property CT for Hida deformations}
\label{completely faithful Hida section}

In this section, we formulate and discuss an analogue of Property CF
and Property CT for Selmer groups of Hida deformations. As before,
$p$ will denote an odd prime. Let $E$ be an elliptic curve over $\Q$
with ordinary reduction at $p$ and assume that $E[p]$ is an
absolutely irreducible $\Gal(\bar{\Q}/\Q)$-representation. By Hida
theory (for instance, see \cite{Hi1,Hi2}; also see \cite{MT,W}),
there exists a commutative complete Noetherian local domain $R$
which is flat over the power series ring $\Zp\ps{X}$ in one
variable, and a free $R$-module $T$ of rank 2 with $T/P \cong T_pE$
for some prime ideal $P$ of $R$. (Note that the freeness of $T$ over
$R$ is guaranteed by our assumption that $E[p]$ is an absolutely
irreducible $\Gal(\bar{\Q}/\Q)$-representation; see \cite{MT} for
details.) \textit{From now on, we will assume that $R = \Zp\ps{X}$
in all our discussion}. The prime ideal $P$ can be expressed as the
ideal generated by $X-a$ for some $a\in p\Zp$. We now mention two
more properties of $T$ which we require to attach an appropriate
Selmer group to the Hida deformation (For more detailed description
of fundamental and important arithmetic properties of the Hida
deformations, we refer readers to \cite{Hi1,Hi2, W}). The first is
that $T$ is unramified outside the set $S$, where $S$ is any finite
set of primes of $F$ which contains the primes above $p$, the
infinite primes, the primes at which $E$ has bad reduction and the
primes that are ramified in $F_{\infty}/F$ (see \cite[Theorem
2.1]{Hi1} or \cite[Theorem 2.2.1]{W}). The second property we will
mention is that there exists an $R$-submodule $T^+$ of $T$ which is
invariant under the action of $\Gal(\bar{\Q}_p/\Qp)$ and such that
both $T^+$ and $T/T^+$ are free $R$-modules of rank one (see
\cite[Theorem 2.2.2]{W} for details; and noting that $T$ is free
over $R$ by our assumption and \cite[Section 2, Corollary 6]{MT}).

Set $A = T\ot_R\Hom_{\cts}(R,\Qp/\Zp)$ and $A^+ =
T^+\ot_R\Hom_{\cts}(R,\Qp/\Zp)$. We note that $E_{p^{\infty}} =
A[P]$. We now define a variant of Selmer group for the Hida
deformation. The notion of such a Selmer group was introduced by
Greenberg in \cite{G89}, where he defined his Selmer group over a
cyclotomic $\Zp$-extension. Since then, his idea has been
adapted to the situation of an admissible $p$-adic Lie extension
$F_{\infty}$ (see \cite[Section 4]{CS12} and \cite[Section 6]{SS}).
Following their footsteps, we define the Selmer group of the Hida
deformation over an admissible $p$-adic Lie extension $F_{\infty}$
of $\Q$ by

\[  S(A/F_{\infty}) = \ker\Big(H^1(G_S(F_{\infty}), A)
\lra \bigoplus_{v\in S} J_v(A, F_{\infty}) \Big), \]
 where
 \[ J_v(A, F_{\infty}) = \begin{cases}
  \prod_{w|v}H^1(F_{\infty,w}, A/A^+),& \mbox{if } v\mbox{ divides }p, \\
      \prod_{w|v}H^1(F_{\infty,w}, A), &
      \mbox{if } v\mbox{ does not divides }p. \end{cases}
 \]
We denote by $X(A/F_{\infty})$ the Pontryagin dual of
$S(A/F_{\infty})$, and consider this dual Selmer group as a
(compact) $\Gal(F_{\infty}/F)$-module for some finite extension $F$
of $\Q$ in $F_{\infty}$, where $F_{\infty}$ is a strongly admissible
$p$-adic Lie extension of $F$.

Let us briefly turn to algebra. For a $R\ps{G}$-module $M$, denote by $M(R)$ the $R$-torsion
submodule of $M$. It is an easy exercise to verify that this is a
$R\ps{G}$-submodule of $M$. We now record the following lemma.

\bl \label{CF R=0}
 Suppose that $M$ is a completely faithful $R\ps{G}$-module. Then
 $M(R)$ is pseudo-null over $R\ps{G}$.
\el

\bpf
 Since $M$ is finitely generated over $R\ps{G}$, we can find $x\in
 R$ such that $M(R) = M[x]$. Since $M$ is completely faithful over
 $R\ps{G}$, this in turn implies that $M(R)$ is pseudo-null over
 $R\ps{G}$.
\epf

Therefore, the above result shows that in considering complete
faithfulness of a $R\ps{G}$-module, we must first cut off its
$R$-torsion submodule. Returning to arithmetic, write $X_R(A/F_{\infty})=
X(A/F_{\infty})/X(A/F_{\infty})(R)$.  We then say that $X_R(A/F_{\infty})$
\textit{satisfies Property CF} if $X_R(A/F_{\infty})$ is either a
pseudo-null $R\ps{G}$-module or a completely faithful
$R\ps{G}$-module.
It will be of interest to
establish a general control theorem relating the validity of
Property CF for the Selmer group of an elliptic curve and the Selmer group of its Hida deformation. At
this point of writing, we are only able to prove such a relation
for a specific class of strongly admissible $p$-adic Lie extensions
(see Theorem \ref{faithful control for Hida}). We now record the
following lemma which is the analogue of Lemma \ref{Coates} and has
a similar proof.

\bl
  Let $F_{\infty}$ be a noncommutative strongly admissible $p$-adic Lie extension of $F$ with
 $G =\Gal(F_{\infty}/F)$. Suppose that $X(A/F_{\infty})$ is finitely generated over $R\ps{H}$. Then the following statements are
 equivalent.
 \begin{enumerate}
 \item[$(a)$] $X(A/F_{\infty})$ is either a
pseudo-null $R\ps{G}$-module or a completely faithful
$R\ps{G}$-module.
  \item[$(b)$] $X_R(A/F_{\infty})$ is either a
pseudo-null $R\ps{G}$-module or a completely faithful
$R\ps{G}$-module.
 \end{enumerate}
\el

To facilitate further discussion,
we recall the
$\M_H(G)$-conjecture for Hida deformations which was first formulated in
\cite{CS12}. We say that $X(A/F_{\infty})$ \textit{satisfies the $\M_H(G)$-Conjecture} if $X_R(A/F_{\infty})$ is a finitely
generated $R\ps{H}$-module.

\medskip
A parallel argument to that in Proposition \ref{faithful Selmer groups}
will establish the following.

\bp \label{faithful Selmer groups hida}
  Let $F_{\infty}$ be a strongly admissible $p$-adic Lie extension of $F$ with
 $G =\Gal(F_{\infty}/F)$. Suppose that the following statements hold.
 \begin{enumerate}
 \item[$(i)$] The $\M_H(G)$ conjecture holds for $X(A/F_{\infty})$.
\item[$(ii)$] $G$ satisfies $\mathbf{(NH)}$.
\end{enumerate}
  Then Property CF holds for $X_R(A/F_{\infty})$. \ep

We now propose an analogue of Property CT for Selmer groups of Hida deformations. As in Section \ref{triviality of central torsion submodules section}, denote by $C$ the center of $G = \Gal(F_{\infty}/F)$, where $F_{\infty}$ is a noncommutative strongly admissible $p$-adic Lie extension of $F$. We then say that $X(A/F_{\infty})$ (or $X_R(A/F_{\infty})$) \textit{satisfies Property CT} if the $R\ps{C}$-torsion submodule of $X_R(A/F_{\infty})$ is zero.
The following analogue of Proposition \ref{CF
CT} has a similar proof which we will omit.

\bp \label{CF CT Hida}
 Let $F_{\infty}$ be a noncommutative strongly admissible $p$-adic Lie extension of $F$. Suppose that all of
  the following statements hold.
 \begin{enumerate}
  \item[$(i)$] The $\M_H(G)$ conjecture holds for
 $X(A/F_{\infty})$.
 \item[$(ii)$] $X(A/F_{\infty})$ has no nonzero pseudo-null $R\ps{G}$-submodule.
\item[$(iii)$] Suppose that either $(a)$ or $(b)$ holds.
 \begin{enumerate}
  \item[$(a)$] $X_R(A/F_{\infty})$ satisfies Property CF.
 \item[$(b)$] The centre $C$ of $G=\Gal(F_{\infty}/F)$ is contained in
 $H$.
\end{enumerate}
\end{enumerate}
 Then $X_R(A/F_{\infty})$ satisfies Property CT.
 \ep

Condition (ii) of the preceding proposition is known to hold in some
situations (see \cite[Theorems 6.10 and 6.12]{SS}). Therefore, one
can combine these observation with Theorem \ref{faithful Selmer
groups hida} to obtain examples of $X_R(A/F_{\infty})$ which
satisfies Property CT. We end this section by mentioning a further interesting corollary which can be viewed as a slight refinement of \cite[Corollary 6.3]{CS12} and
\cite[Proposition 8.2]{SS}.

\bc \label{CF CT2 Hida}
  Retain the assumptions of Proposition \ref{CF CT Hida}. Suppose further
  that $X(A/F_{\infty})$ is finitely generated over
 $R\ps{H}$ with positive $R\ps{H}$-rank.
 Then the $R\ps{C}$-torsion submodule of $X(A/F_{\infty})$ is zero.  \ec

\section{Control theorems from Property CF to Property CT} \label{control CT section}

In this section, we will prove our main result. Namely, we will prove a control theorem which allows us to
conclude that the dual Selmer group satisfies Property CT if we know
that the dual Selmer group satisfies Property CF in a smaller field.
We also establish an analogous result for Hida deformation which says
that if the dual Selmer group of an elliptic curve satisfies Property CF,
then the Selmer group of its Hida deformation satisfies Property CT. As a start, we begin with an algebraic
result.

\bl \label{control alg prop} Let $G$ be a nonabelian pro-$p$
$p$-adic Lie group with no $p$-torsion. Suppose that $G$ contains
two normal closed subgroups $N$ and $H$ such that $N\subseteq H$,
$G/H\cong\Zp$, $N\cong \Zp$ and $G/N$ is a non-abelian pro-$p$
$p$-adic Lie group without $p$-torsion. Let $M$ be a
$\Zp\ps{G}$-module which satisfies all the following properties.
\begin{enumerate}
\item[$(i)$] $M$ is finitely generated over $\Zp\ps{H}$.
\item[$(ii)$] $M$ has
no nonzero pseudo-null $\Zp\ps{G}$-submodule.
\item[$(iii)$] $H_1(N,M)=0$.
\item[$(iv)$] $M_N$ is either pseudo-null or
completely faithful over $\Zp\ps{G/N}$.
\end{enumerate}
 Then the $\Zp\ps{C}$-torsion
submodule of $M$ is zero, where $C$ denotes the center of $G$. \el

\bpf
   We first suppose that $M_N$ is a pseudo-null
  $\Zp\ps{G/N}$-module. Since $M$ is finitely generated over $\Zp\ps{H}$, its $N$-invariant
  $M_N$ is certainly finitely generated over $\Zp\ps{H/N}$. It then follows from Theorem \ref{pseudo-null
  torsion} that $M_N$ is a finitely generated torsion
  $\Zp\ps{H/N}$-module.  By Lemma \ref{relative rank}, we have
  \[ 0 = \rank_{\Zp\ps{H/N}} M_{N} =  \rank_{\Zp\ps{H}} M  +
 \rank_{\Zp\ps{H/N}} H_1(N,M)  \]
 which in turn implies that $\rank_{\Zp\ps{H}}M=0$.
 Hence $M$ is a finitely generated torsion
  $\Zp\ps{H}$-module, or equivalently, a pseudo-null
  $\Zp\ps{G}$-module. But since $M$ has
no nonzero pseudo-null $\Zp\ps{G}$-submodule, we have
$M=0$, and in particular, the $\Zp\ps{C}$-torsion submodule of $M$ is
zero.

We now suppose that $M_N$ is a completely faithful
$\Zp\ps{G/N}$-module. Write $W$ for the $\Zp\ps{C}$-torsion submodule of $M$. Since $C$ is the center of $G$ and $M$ is finitely generated over $\Zp\ps{G}$, we have $W \subseteq M[z]$ for some $z\in \Zp\ps{C}$. The reverse inclusion is evident from the definition of $W$ being the $\Zp\ps{C}$-torsion submodule of $M$. As $z$ is a central element of $\Zp\ps{G}$, the multiplication by $z$-map on $M$ is a $\Zp\ps{G}$-homomorphism. This homomorphism in turn gives an identification $zM \cong M/M[z] = M/W$, where the last equality follows from the above discussion. It then follows from this that $M/W$ can be viewed as a $\Zp\ps{G}$-submodule of $M$. Since $H_1(N,M)=0$
by hypothesis and $H_1(N,-)$ is left exact, we have $H_1(N,M/W)=0$. The left exactness of $H_1(N,-)$ also yields $H_1(N, W)=0$ which is equivalent to saying that $W[\ga_N-1] = 0$, where $\ga_N$ is a topological generator of $N$. Note that the augmentation kernel $I_N$ of $\Zp\ps{G}\lra \Zp\ps{G/N}$ is generated by $\ga_N-1$. Since $I_N=\Zp\ps{G}(\ga_N-1)=(\ga_N-1)\Zp\ps{G}$, we may therefore write $z = z'(\ga_N-1)^m$ for some $z'\notin I_N$ and nonnegative integer $m$. By virtue of $W[\ga_N-1] = 0$, one checks easily that $W= M[z] = M[z']$. Hence we may assume that $W= M[z]$ with $z\notin I_N$ which we do. Therefore, the image of $z$ in $\Zp\ps{G/N}$, which we denote by $\bar{z}$, is nonzero. Now, applying $N$-invariant to the short exact sequence
  \[ 0\lra W\lra M \lra M/W \lra 0,\]
  we obtain a short exact sequence
  \[ 0\lra W_N\lra M_N \lra (M/W)_N \lra 0,\]
  where the injectivity of the leftmost map is a consequence of the above observation that $H_1(N,M/W)=0$. Clearly, $\bar{z}$ is a global annihilator of $W_N$. But since $M_N$ is assumed to be a completely faithful $\Zp\ps{G/N}$-module, this forces
  $W_N$ to be a pseudo-null $\Zp\ps{G/N}$-module. On the other hand,
  since $M$, and hence $W$, is finitely generated over $\Zp\ps{H}$, we have that $W_N$ is finitely generated over $\Zp\ps{H/N}$. It then follows from Theorem \ref{pseudo-null torsion}
  that $W_N$ is a torsion $\Zp\ps{H/N}$-module. By Lemma \ref{relative rank}, this implies that
  $W$ is a torsion $\Zp\ps{H}$-module. Now, by another application of Theorem \ref{pseudo-null torsion}, we have that $W$ is pseudo-null over $\Zp\ps{G}$. Since $M$ has no nonzero pseudo-null $\Zp\ps{G}$-submodule by the hypothesis of the lemma, we conclude that $W =0$, completing the proof of the lemma. \epf

We can now prove the first of our main theorem. Before doing so, we introduce further
notation. Let $F_{\infty}$ be a strongly admissible $p$-adic Lie
extension of $F$ and $E$ an elliptic curve over $F$.  For the remainder of the paper,
$S$ will denote a fixed finite set of primes of $F$ which contains the primes above $p$, the
infinite primes, the primes at which $E$ has bad reduction and the
primes that ramify in $F_{\infty}/F$. We then denote by $S_p$ the subset
of $S$ consisting of primes of $F$ above $p$ and by
$S_2(F_{\infty}/F)$ the subset of $S$ consisting of primes of $F$
whose decomposition group of $\Gal(F_{\infty}/F)$
  at $w$ has dimension $\geq 2$, where $w$ is a prime of
  $F_{\infty}$ above $v$.

\bt \label{faithful trivial torsion}
 Let $E$ be an elliptic curve over $F$
with either good ordinary reduction or multiplicative reduction at
every prime of $F$ above $p$.
 Let $F_{\infty}$ be a strongly admissible $p$-adic Lie extension of $F$ with
 $G= \Gal(F_{\infty}/F)$.
  Suppose that the following statements hold.
  \begin{enumerate}
\item[$(i)$] $N$ is a closed normal subgroup of $G$ which is contained
in $H$ and $N\cong \Zp$.
\item[$(ii)$] $G/N$ is a non-abelian pro-$p$ $p$-adic Lie group without
$p$-torsion. $($In particular, the dimension of the $p$-adic Lie
group $G/N$ is necessarily $\geq 2$.$)$
\item[$(iii)$]  Suppose that the $\M_H(G)$ conjecture holds for
$X(E/F_{\infty})$ and that $X(E/F_{\infty})$ has no nonzero
pseudo-null $\Zp\ps{G}$-submodules.
\item[$(iv)$]  $X_f(E/L_{\infty})$ satisfies Property CF. Here $L_{\infty} : = F_{\infty}^N$.
\item[$(v)$]  Suppose that either $(a)$ or $(b)$ holds.
 \begin{enumerate}
\item[$(a)$] $S = S_2(L_{\infty}/F)$.
\item[$(b)$] $F$ contains a primitive $p$-root of unity and $S_p\subseteq S_2(L_{\infty}/F)$. Also, suppose that for every
 $v\in S\setminus S_2(L_{\infty}/F)$, either the
decomposition group of $\Gal(F_{\infty}/F)$ at $v$ has dimension 1,
or $E$ has non-split multiplicative reduction or additive reduction
at $v$. Furthermore, assume that the elliptic curve $E$ has no
additive reduction when $p = 3$.
\end{enumerate}
\end{enumerate}
 Then $X_f(E/F_{\infty})$ satisfies Property CT. \et

\br
 Of course, if the center $C$ is contained in $H$, then the
 conclusion of the theorem is immediate from Proposition \ref{CF
 CT}. The point of the above theorem is that we do not assume that
 the center $C$ is contained in $H$.
\er

\bpf[Proof of Theorem \ref{faithful trivial torsion}]
  It suffices to verify that all the hypothesis of
Lemma \ref{control alg prop} are satisfied (taking
$M=X_f(E/F_{\infty})$). The first assertion of assumption (iii) guarantees us that
$X_f(E/F_{\infty})$ is finitely generated over $\Zp\ps{H}$. Combining the second assertion of assumption (iii) with \cite[Lemma 4.2]{Su}, we have that
$X_f(E/F_{\infty})$ has no nonzero pseudo-null
$\Zp\ps{G}$-submodules. By the argument in \cite[Remark 6.2]{LimCF},
we have $H_1(N, X(E/F_{\infty}))=0$. Since $X_f(E/F_{\infty}) =
p^mX(E/F_{\infty})$ for a large enough $m$, we also have $H_1(N,
X_f(E/F_{\infty}))=0$ by the left exactness of $H_1(N,-)$.
Therefore, it remains to show that $X_f(E/F_{\infty})_N$ is either
pseudo-null or completely faithful over $\Zp\ps{G/N}$. By a standard well-known
argument (for instance, see \cite[Lemma 2.4]{CS12}), one can show that the map
  \[ \al: X(E/F_{\infty})_N \lra X(E/L_{\infty}) \]
has kernel which is finitely generated over $\Zp\ps{H/N}$, and
cokernel which is finitely generated over $\Zp$. We further claim
that under the hypothesis of the theorem, $\ker \al$ is a finitely generated torsion
$\Zp\ps{H/N}$-module, and therefore, is pseudo-null over
$\Zp\ps{G/N}$. To prove our claim, we first note that since the Pontryagin dual of the map
\[ a:S(E/L_{\infty})\lra S(E/F_{\infty})^N\]
is precisely the map
\[ \al:X(E/F_{\infty})_N\lra X(E/L_{\infty}),\]
 we are reduced to showing
that $\coker a$ is a cofinitely generated torsion
$\Zp\ps{H/N}$-module. Now consider the following commutative diagram
\[
 \entrymodifiers={!! <0pt, .8ex>+} \SelectTips{eu}{}\xymatrix{
    0 \ar[r]^{} & S(E/L_{\infty}) \ar[d]_{a} \ar[r] &
    H^1(G_S(L_{\infty}), E_{p^{\infty}})
    \ar[d]_{b}
    \ar[r] & \bigoplus_{v\in S}J_v(E/L_{\infty}) \ar[d]_{c=\oplus c_v}\ar[r] &0\\
    0 \ar[r]^{} & S(E/F_{\infty})^N\ar[r]^{}
    & H^1(G_S(F_{\infty}),  E_{p^{\infty}})^N \ar[r] & \
    \bigoplus_{v\in S}J_v(E/F_{\infty})^N  & } \]
with exact rows. (The surjectivity of the last map in the top exact sequence can be shown by a similar argument in the proof of Lemma \ref{observe}.) Since $N\cong\Zp$, the map $b$ is surjective. Thus, it suffices via the snake lemma to show that $\ker c$ is a cofinitely generated torsion $\Zp\ps{H/N}$-module. Since $S$ is finite, we are reduced to showing that $\ker c_v$ is a cofinitely generated torsion $\Zp\ps{H/N}$-module for every $v\in S$. To see this, we first observe that one can apply a similar
argument in the spirit of the proof of
\cite[Lemma 8.7]{SS} to show that $\ker c_v$ is a cofinitely generated
torsion $\Zp\ps{H/N}$-module for $v\in S_2$. It therefore remains to
consider the case when $v\in S\setminus S_2$ in assumption (iv)(b).
Write $c_v = \oplus c_w$, where the sum runs over the primes of
$L_{\infty}$ above $v$. It is then an easy exercise to verify that
$\ker c_w = H_1(N_u, E(F_{\infty,u})_{p^{\infty}})$, where $u$ is a
prime of $F_{\infty}$ above $w$ and $N_u$ is the decomposition group
of $N= \Gal(F_{\infty}/L_{\infty})$ at $u$. Now suppose that the
decomposition group of $\Gal(F_{\infty}/F)$ at $v$ has dimension 1.
Since $v$ does not divide $p$ (because $S_p\subseteq
S_2(L_{\infty}/F)$) and $F_{\infty}$ contains $F^{\cyc}$, we have
$F_{\infty,u}= F^{\mathrm{ur},p}_v$, where $F^{\mathrm{ur},p}_v$ is
the maximal unramified pro-$p$ extension of $F_v$. Similarly, we
have $L_{\infty,w}= F^{\mathrm{ur},p}_v$. Hence this implies that
$L_{\infty,w}=F_{\infty,w}$, or in other words,
$N_u =0$. In particular, we have $\ker c_v=0$ when the decomposition group of $\Gal(L_{\infty}/F)$
at $v$ has dimension $1$. Now suppose that $E$ has non split
multiplicative reduction or additive reduction at $v$. Since $p$ is odd (and $p\geq 5$ if $E$ has additive reduction), $E$ also has non split multiplicative reduction or additive reduction at every primes $x$ of
$F^{\cyc}$ above $v$. It then follows from \cite[Proposition 5.1(iii)]{HM} that $E(F^{\cyc}_x)_{p^{\infty}} =0$. Since  $F_{\infty, w}/F^{\cyc}_x$ is a pro-$p$ extension, we have $E(F_{\infty, w})_{p^{\infty}}
=0$ and thus one also has $\ker c_v=0$ in this case. In conclusion, we have that $\ker c$ is a cofinitely generated torsion $\Zp\ps{H/N}$-module. This completes the verification of our claim that $\ker \al$ is a finitely generated torsion
$\Zp\ps{H/N}$-module.

Now consider the following commutative diagram
 \[ \entrymodifiers={!! <0pt, .8ex>+}
\SelectTips{eu}{}\xymatrix{
    0 \ar[r]^{} & X(E/F_{\infty})(p)_N \ar[d]_{\al'} \ar[r] &
    X(E/F_{\infty})_N
    \ar[d]_{\al}
    \ar[r] & X_f(E/F_{\infty})_N \ar[d]_{\al''} \ar[r] & 0 \\
    0 \ar[r]^{} & X(E/L_{\infty})(p) \ar[r]^{} & X(E/L_{\infty})\ar[r] & \
    X_f(E/L_{\infty}) \ar[r] &  0  } \]
with exact rows (Here the top exact row follows from the above
mentioned fact that $H_1(N, X_f(E/F_{\infty}))=0$). Since $\coker
\al$ is finitely generated over $\Zp$ and the dimension of $G/N$ is $\geq 2$, we have that $\coker\al$ is pseudo-null over $\Zp\ps{G/N}$. It then follows from the above commutative diagram that $\coker \al''$
is pseudo-null over $\Zp\ps{G/N}$. We
shall now show that $\ker \al''$ is pseudo-null over $\Zp\ps{G/N}$.
By assumption (iii) of the theorem, we have that $X_f(E/F_{\infty})_N$ is
finitely generated over $\Zp\ps{H/N}$ and this in turn implies that
$\ker \al''$ is finitely generated over $\Zp\ps{H/N}$. Combining
this with the observation that $\coker \al$ is finitely
generated over $\Zp$, we have that $\coker \al'$ is finitely
generated over $\Zp\ps{H/N}$. Since $\coker \al'$ is $p$-primary, it
follows from Lemma \ref{mu lemma}(b) and (c) that $\coker\al'$ is
pseudo-null over $\Zp\ps{G/N}$. Now, combining this with the above
proven claim that $\ker \al$ is pseudo-null over $\Zp\ps{G/N}$,
we have that $\ker \al''$ is pseudo-null over $\Zp\ps{G/N}$. Hence we conclude that the map
\[ \al'':X_f(E/F_{\infty})_N\lra X_f(E/L_{\infty})\]
has kernel and cokernel that are pseudo-null over $\Zp\ps{G/N}$.
Since $X_f(E/L_{\infty})$ satisfies Property CF, it then follows from an appplication of
Lemma \ref{cf compare} that $X_f(E/F_{\infty})_N$ is either
pseudo-null or completely faithful over $\Zp\ps{G/N}$. The proof of Theorem \ref{completely faithful
control} is now complete.
 \epf

We now consider a variant of Theorem \ref{faithful trivial torsion}
for an elliptic curve and its Hida deformation. Recall that $A$ is
the $R$-cofree Galois module attached to the Hida deformation as
defined in Section \ref{completely faithful Hida section}, where
$R=\Zp\ps{X}$, and has the property that $A[P] = E_{p^{\infty}}$ for
some prime ideal $P$ of $R$ generated by $X-a$ with $a\in p\Zp$. As
before, we denote by $X(A/F_{\infty})$ the dual Selmer group of the
Hida deformation.

\bt \label{faithful control for Hida CT}
 Let $F_{\infty}$ be a strongly admissible $p$-adic Lie extension of $F$
 with Galois group
 $G$. Suppose that the following statements hold.
     \begin{enumerate}
     \item[$(i)$] $X(A/F_{\infty})$ satisfies the $\M_H(G)$ conjecture.
\item[$(ii)$]   $X(A/F_{\infty})$ has no nonzero pseudo-null
$\Zp\ps{G}$-submodules.
\item[$(iii)$] For every $v\in S$, the decomposition group of $\Gal(F_{\infty}/F)$
  at $v$ has dimension $\geq 2$.
\item[$(iv)$] $X_f(E/F_{\infty})$ satisfies Property CF.
\end{enumerate}
 Then  $X_R(A/F_{\infty})$ satisfies Property CT.
\et

We first record a preliminary lemma. Recall that a polynomial $X^n +
c_{n-1}X^{n-1} + \cdots +c_0$ in $\Zp[X]$ is said to be a
\textit{Weierstrass polynomial} if $p$ divides $c_i$ for every
$0\leq i \leq n-1$.

\bl
 Let $M$ be a finitely generated $R\ps{G}$-module with the property
 that $M = M(R)$. If $P$ is a prime ideal of $R$ generated by $X-a$ for some $a\in p\Zp$ with the
 property that $M[X-a] = 0$. Then $M/PM$ is a finitely generated
 $p$-primary $\Zp\ps{G}$-module.
\el

\bpf
 The finite generation of $M$ over $R\ps{G}$ allows us to find a Weierstrass
polynomial $f(X)$ such that $f(X)M =0$.  Now write $f(X) =
(X-a)^mg(X)$, where $g(a)\neq 0$. Since $M[X-a] = 0$, one easily
checks that $g(X) M = 0$. Clearly, $g(a)$ is nonzero and annihilates
$M/PM$. Since $a\in p\Zp$ and $g(X)$ is a Weierstrass
polynomial, we have $g(a)\in p\Zp$. This shows that $M/PM$ is
$p$-primary.   \epf

We can now give the proof of Theorem \ref{faithful control for Hida
CT}.

\bpf[Proof of Theorem \ref{faithful control for Hida CT}]
 The proof is essentially similar to that in Theorem \ref{faithful trivial torsion},
 where we choose an identification $R\ps{G} \cong \Zp\ps{N\times G}$ such that
 $M/PM = M_N$ for every $R\ps{G}$-module $M$. Here $N\cong \Zp$.
 The only thing which perhaps requires additional attention
 is to establish the following commutative diagram
 \[ \entrymodifiers={!! <0pt, .8ex>+}
\SelectTips{eu}{}\xymatrix{
    0 \ar[r]^{} & X(A/F_{\infty})(R)/P\ar[d]_{\al'} \ar[r] &
    X(A/F_{\infty})/P
    \ar[d]_{\al}
    \ar[r] & X_R(A/F_{\infty})/P \ar[d]_{\al''} \ar[r] & 0 \\
    0 \ar[r]^{} & X(E/F_{\infty})(p) \ar[r]^{} & X(E/F_{\infty})\ar[r] & \
    X_f(E/F_{\infty}) \ar[r] &  0  } \]
with exact rows. By the preceding lemma, $X(A/F_{\infty})(R)/P$ is
$p$-primary (noting that $X(A/F_{\infty})(R)[X-a] =
H_1(N,X(A/F_{\infty})(R)) =0$ by a similar argument to that in
Theorem \ref{faithful trivial torsion}). This in turn implies that
$X(A/F_{\infty})(R)/P$ is sent into $X(E/F_{\infty})(p)$ under
$\al$, and thus inducing the required maps $\al'$ and $\al''$. The
remainder of the proof proceeds as in  Theorem \ref{faithful trivial
torsion}. \epf

\section{Control theorems for Property CF} \label{control section}

In \cite[Section 6]{LimCF}, the author has established control
theorems for the faithfulness of dual Selmer groups. It will be
of interest to obtain analogous results for Property CF. In this
section, we will establish such a control theorem for Property CF over
certain class of strongly admissible $p$-adic Lie extensions.
We first prove the control theorem which concerns
extensions of strongly admissible $p$-adic Lie extensions.
\textbf{We will assume that $p\geq 5$ throughout the section.}

As a start, we record the following theorem of Csige \cite[Theorem
2.4.1]{Cs} which generalizes a previous result of Ardakov
\cite[Theorem 1.3]{A}.

\bt[Ardakov, Csige] \label{Csige} Let $p\geq 5$. Given that $G \cong
\Zp^r\times H_0 \times \Ga$. Here $r\geq 0$, and $H_0$ is a torsion-free
compact $p$-adic analytic group, whose Lie algebra is split
semisimple over $\Qp$. Let $M$ be $\Zp\ps{G}$-module which has no
nonzero pseudo-null $\Zp\ps{G}$-submodules. Write $C=  \Zp^r\times
\Ga$. Then $M$ is a completely faithful $\Zp\ps{G}$-module if and
only if the $\Zp\ps{C}$-torsion submodule of $M$ is trivial. \et

We can now state and prove our control theorems.

\bt \label{completely faithful control}
 Let $p\geq 5$. Let $E$ be an elliptic curve over $F$
with either good ordinary reduction or multiplicative reduction at
every prime of $F$ above $p$.
 Let $F_{\infty}$ be a strongly admissible $p$-adic Lie extension of $F$ with
 $G= \Gal(F_{\infty}/F)$.
  Suppose that the following statements hold.
  \begin{enumerate}
\item[$(i)$] $X(E/F_{\infty})$ satisfies the $\M_H(G)$ conjecture.
\item[$(ii)$] $G= H\times \Ga$, where $H\cong N\times \Zp^r\times
H_0$. Here $r\geq 0$, $H_0$ is a torsion-free compact $p$-adic
analytic group, whose Lie algebra is split semisimple over $\Qp$ and
$N\cong \Zp$.
\item[$(iii)$] $X(E/F_{\infty})$ has no nonzero pseudo-null
$\Zp\ps{G}$-submodules.
\item[$(iv)$] Set $L_{\infty} : = F_{\infty}^N$. For every prime $v$ of $F$ above $p$, the decomposition group of $\Gal(L_{\infty}/F)$
  at $v$ has dimension $\geq 2$.
\item[$(v)$]  $X_f(E/L_{\infty})$ satisfies Property CF.
\end{enumerate}
  Then $X_f(E/F_{\infty})$ satisfies Property CF. \et

\bpf
 By Theorem \ref{Csige}, it suffices to show that the $\Zp\ps{C}$-torsion submodule of
 $X_f(E/F_{\infty})$ is trivial. By Lemma \ref{control alg
 prop}, we are then reduced to showing that $X_f(E/F_{\infty})_N$ is either
pseudo-null or completely faithful over $\Zp\ps{G/N}$. The proof is
essentially similar to that in Theorem \ref{faithful trivial
torsion}. As the decomposition group condition in this theorem
differs slightly from those in Theorem \ref{faithful trivial
torsion}, the only thing which perhaps requires additional attention
is to verify that $\ker c$ (see Theorem \ref{faithful trivial torsion} for notation) is a cofinitely generated torsion $\Zp\ps{H/N}$-module. For each $v\in S$ such that the decomposition group of $\Gal(L_{\infty}/F)$ at $v$ has dimension $\geq 2$, we can apply a similar argument in the spirit of the proof of \cite[Lemma 8.7]{SS} to show that $\ker c_v$ is a cofinitely generated
torsion $\Zp\ps{H/N}$-module. Let $v\in S$ such that the
decomposition group of $\Gal(L_{\infty}/F)$ at $v$ has dimension
$1$. Then $v$ does not divide $p$ by assumption (iv). Let $w$ be a
prime of $F_{\infty}$ above $v$. By abuse of notation, we also
denote by $w$ the prime of $L_{\infty}$ below $w$. As the
decomposition group of $\Gal(L_{\infty}/F)$ at $v$ has dimension
$1$, we have $L_{\infty,w} = F^{\cyc}_w = F^{\mathrm{ur},p}_v$,
where $F^{\mathrm{ur},p}_v$ is the maximal unramified pro-$p$ extension
of $F_v$. Since $F_{\infty}$ is obtained from $L_{\infty}$ by adjoining a
$\Zp$-extension of $F$, it follows that $F_{\infty}/L_{\infty}$ is
also unramified outside $p$. Since $L_{\infty,w}$ is the maximal unramified pro-$p$ extension
of $F_v$, we have $L_{\infty,w}=F_{\infty,w}$, or in other words, we have $N_w =0$, where $N_w$ is the decomposition group of $N= \Gal(F_{\infty}/L_{\infty})$ at $w$. Hence we have $\ker c_v=0$ if the decomposition group of $\Gal(L_{\infty}/F)$ at $v$ has dimension $1$. In conclusion, we have shown that $\ker c$ is a cofinitely generated torsion $\Zp\ps{H/N}$-module. The remainder of the
proof proceeds as in  Theorem \ref{faithful trivial torsion}. \epf

The following is an analogue of Theorem \ref{completely faithful
control} for an elliptic curve and its Hida deformations which has a
similar proof that we will omit.

\bt \label{faithful control for Hida}
 Let $p\geq 5$. Let $F_{\infty}$ be a strongly admissible $p$-adic Lie extension of $F$
 with Galois group
 $G$. Suppose that the following statements hold.
     \begin{enumerate}
\item[$(i)$] The $\M_H(G)$ conjecture holds for $X(A/F_{\infty})$.
\item[$(ii)$] $G= H\times \Ga$, where $H\cong \Zp^r\times
H_0$. Here $r\geq 0$, and $H_0$ is a torsion-free compact $p$-adic
analytic group, whose Lie algebra is split semisimple over $\Qp$.
\item[$(iii)$] $X(A/F_{\infty})$ has no nonzero pseudo-null
$R\ps{G}$-submodules.
\item[$(iv)$] For every $v\in S$, the decomposition group of $\Gal(F_{\infty}/F)$
  at $v$ has dimension $\geq 2$.
\item[$(v)$] $X_f(E/F_{\infty})$ satisfies Property CF.
\end{enumerate}
 Then $X_R(A/F_{\infty})$ satisfies Property CF.
\et

\section{Examples} \label{examples section}

In this section, we will give some examples to illustrate the
results in this paper.

\medskip
(a) Let $E$ be the elliptic curve $11a2$ of Cremona's table which is
given by
 \[ y^2 + y = x^3 -x. \]
Take $p = 5$, $F=\Q(\mu_5)$ and $L_{\infty} = \Q(\mu_{5^{\infty}},
11^{5^{-\infty}})$. Let $F_{\infty}$ be a strongly admissible
$5$-adic Lie extension of $F$ that contains $L_{\infty}$ and that
the group $N=\Gal(F_{\infty}/L_{\infty})$ satisfies the conditions
in Theorem \ref{faithful Selmer groups}. As observed in \cite{HV,
LimCF}, $X(E/F_{\infty})$, and hence $X_f(E/F_{\infty})$, is
finitely generated over $\Z_5\ps{H}$ with positive
$\Z_5\ps{H}$-rank. By \cite[Proposition 5.6]{CFKSV}, the
$\M_H(G)$-conjecture holds for $X(E'/F_{\infty})$, where $E'$ is
either $11a1$ or $11a3$. From the proof of Proposition \ref{faithful
isogeny}, we also see that $X_f(E'/F_{\infty})$ has the same
$\Z_5\ps{H}$-rank as $X(E/F_{\infty})$. Proposition \ref{faithful
Selmer groups} then tells us that $X_f(A/F_{\infty})$ is completely
faithful over $\Z_5\ps{G}$, and in particular satisfies Property CF
for $A=11a1$, $11a2$ and $11a3$.

Furthermore, it follows from an application of Proposition \ref{CF
CT2} that the $\Z_5\ps{C}$-torsion submodules of $X_f(A/F_{\infty})$
are zero, where $C$ is the center of $G = \Gal(F_{\infty}/F)$. Some
interesting examples of strongly admissible $5$-adic extensions
$F_{\infty}$ of $F$ that one can take are
\[M_{\infty}(\mu_{5^{\infty}}, 2^{5^{-\infty}},11^{5^{-\infty}}),\quad
M_{\infty}(\mu_{5^{\infty}}, 3^{5^{-\infty}}, 11^{5^{-\infty}}),
\quad M_{\infty}(\mu_{5^{\infty}}, 2^{5^{-\infty}}, 3^{5^{-\infty}},
11^{5^{-\infty}}),\]
\[M_{\infty}(E_{5^{\infty}}, 2^{5^{-\infty}},11^{5^{-\infty}}),\quad
M_{\infty}(E_{5^{\infty}}, 3^{5^{-\infty}}, 11^{5^{-\infty}}), \quad
M_{\infty}(E_{5^{\infty}}, 2^{5^{-\infty}}, 3^{5^{-\infty}},
11^{5^{-\infty}}),\] where $M_{\infty}$ is any $\Z_5^r$-extension of
$F$. We should mention that we do not have an answer to the validity of Property CF for the
admissible $5$-adic extensions in the second line.

\medskip
(b) Let $p = 5$. Let $E$ be the elliptic curve $21a4$ of Cremona's
tables which is given by
 \[ y^2 + xy = x^3 + x.\]
Let $p= 5$ and $F=\Q(\mu_5)$. As discussed in \cite[Section 7]{BZ},
$X(E/F)$ is finite as suggested by its $p$-adic $L$-function which
in turn implies that $X(E/F^{\cyc})$ has trivial
$\mu_{\Ga}$-invariant and $\la$-invariant. \textit{We will assume
this latter property throughout our discussion here}.  Let $A$ be an
elliptic curve $66c1$ of Cremona's tables given by
 \[A: y^2 + xy = x^3 - 45x  - 81. \]
Set $L_{\infty} = F(A_{5^{\infty}})$. Since $A$ has no complex
multiplication and $A(\Q)_{\mathrm{tor}}\cong \Z/2\times \Z/5$ (for
instance, see \cite{Cr}), it follows that $L_{\infty}$ is a
noncommutative pro-$p$ extension of $F$. Since $p\geq 5$,
$L_{\infty}$ is in fact a strongly admissible $p$-adic Lie extension
of $F$.

Let $P_0, P_1, P_2$ be defined as in \cite[Definition 1.2]{BZ}. Then
in this situation, we have $P_0 = \{2,3,11\}$. Since $E$ has split
multiplicative reduction at 3 and 3 does not split in $F/\Q$, we
have $P_1 =\{3\}$. We shall now show that $P_2$ is empty. By the discussion in \cite[Section 7]{BZ}, we
have that 2 does not lie in $P_2$.
Therefore, it remains to show that 11 does not lie in $P_2$. Since
the prime 11 splits completely over $F/\Q$, we have $F_u = \Q_{11}$
for every prime $u$ of $F$ above 11. Since $E$ has good reduction at
5, we have an injection $E(\Q_{11})[5] \hookrightarrow
E(\mathbb{F}_{11})$ (see \cite[Chap. VII., Proposition 3.1(b)]{Si}).
By inspection, one has $|E(\mathbb{F}_{11})| = 8$. In particular, this 
implies that $E(\Q_{11})[5] = 0$, and thus all the primes of $F$ above
11 do not lie in $P_2$. (Alternatively, to see that $E(\Q_{11})[5] =
0$, one may give a proof by contradiction. Suppose that
$E(\Q_{11})[5]\neq 0$. Then since $|E(\Q)|=4$, it follows from the
injections $E(\Q)\hookrightarrow E(\mathbb{F}_{11})$ and
$E(\Q_{11})[5]\hookrightarrow E(\mathbb{F}_{11})$ that 20 divides
$|E(\mathbb{F}_{11})|$. But by the Hasse bound (cf. \cite[Chap. V.,
Theorem 1.1]{Si}), we have
\[|E(\mathbb{F}_{11})| \leq 11 + 1 + 2\sqrt{11} = 18.633... <20\]
and this gives the required contradiction.) We may now apply
\cite[Corollary 6.3]{BZ} to conclude that $X(E/L_{\infty})$ is
completely faithfully over $\Z_5\ps{\Gal(L_{\infty}/F)}$, noting that
$X(E/L_{\infty})$ is finitely generated over
$\Z_5\ps{\Gal(L_{\infty}/F^{\cyc})}$).

Let $F_{\infty}$ be the compositum of $L_{\infty}$ and $M_{\infty}$,
where $M_{\infty}$ is any $\Z_5^r$-extension of any finite extension
$L$ of $F$ contained in $L_{\infty}$. All the hypothesis in Theorem
\ref{completely faithful control} can be verified easily. We will
just make the remark that assumption (iii) follows from \cite[Remark
2.6 and Corollary 6.1(on p.\ 1243)]{BZ}, and assumption (iv) follows
from the fact that $A$ has good ordinary reduction at 5 and
\cite[Lemma 2.8(ii)]{C}. Therefore, we have that $X(E/F_{\infty})$
is a completely faithful $\Z_5\ps{\Gal(F_{\infty}/F)}$-module by an
iterative application of Theorem \ref{completely faithful control}.
One can then apply Proposition \ref{faithful isogeny} to show that
$X_f(E'/F_{\infty})$ is also completely faithful over $\Z_5\ps{G}$
for $E' =  21a1, 21a2, 21a3, 21a5, 21a6$. Therefore, Property CF
holds for the Selmer groups of $E$ and its isogenous elliptic curves
over such $F_{\infty}$. By Propositions \ref{trivial torsion
isogeny} and \ref{CF CT}, Property CT also holds for the Selmer
groups of $E$ and its isogenous elliptic curves over such
$F_{\infty}$.

\begin{ack}
    The author would like to thank Romyar Sharifi for his comments and
    for pointing out some grammatical issues in the paper. The author is supported by the
     National Natural Science Foundation of China under the Research
Fund for International Young Scientists
     (Grant No: 11550110172).
    \end{ack}

\end{document}